# Linear Hypergraph List Edge Coloring – Generalizations of the EFL Conjecture to List Coloring

Vance Faber[1, a *]

[1]Center for Computing Sciences, Bowie, Maryland

[a]vance.faber@gmail.com

Revision: January 3, 2017

**Abstract.** Motivated by the Erdős-Faber-Lovász (EFL) conjecture for hypergraphs, we consider the list edge coloring of linear hypergraphs. We discuss several conjectures for list edge coloring linear hypergraphs that generalize both EFL and Vizing's theorem for graphs. For example, we conjecture that in a linear hypergraph of rank 3, the list edge chromatic number is at most 2 times the maximum degree plus 1. We show that for sufficiently large fixed rank and sufficiently large degree, the conjectures are true.

**Generalizations**

**Preliminaries.** This paper is an adaption of [1] to list coloring. Before we can discuss extensions to EFL, we need to give a short list of the concepts involved.

   **Notation**. Let $H = (V, E)$ be a *hypergraph* (see, for example [2]): a set of subsets $E$ of the set $V$. We call the elements of $V$ the *vertices* and the elements of $E$ the *edges*. We often write $n = |V|$ and $m = |E|$. The *degree* of a vertex $x$ is the number of edges $d(x)$ which include it. We let the minimum degree be $\delta$ and the maximum degree be $\Delta$. If all vertices have the same degree, we say the hypergraph is *regular*. The *rank* of an edge $e$ is the cardinality $r(e)$ of $e$. We let the minimum rank be $\rho$ and the maximum rank be $\mathrm{P}$. If all edges have the same rank, we say the hypergraph is *uniform*. If $r(e) = 2$ for every edge then $H$ is a *graph*. If the intersection of any two edges has at most one vertex, we call the hypergraph *linear*.

   **Incidence matrix formulation**. An equivalent formulation for a hypergraph is to consider $H$ to be the incidence matrix of the hypergaph. In this case, $H$ is an $n \times m$ matrix: a row of $H$ is the transpose of the characteristic vector of a vertex and a column of $H$ is the characteristic vector of a edge. We use these two formulations interchangeably. It is often easier to understand a fact in one formulation or the other. For example, a fundamental theorem for a hypergraph is that the sum of the ranks is equal to the sum of the degrees. This is trivial to see in the matrix formulation because both sides of the equality are clearly equal to the number of non-zero entries in the matrix $H$. In this formulation, an edge $e$ is a column vector and a vertex $x$ is a row vector. Two vertices $x$ and $y$ are independent if and only if they are orthogonal, that is, the inner product $xy^* = (x^*, y^*) = 0$. Two edges are independent if and only if they are orthogonal, that is $e^* f = (e, f) = 0$.

   **Duality**. If $H$ is a hypergraph, so is the transpose $H^*$ called the *dual* hypergraph. The theorems we need are often stated in a form that applies most naturally to $H^*$ and we have to translate them to $H$ to apply them. Clearly the edges of $H^*$ are the vertices of $H$ and vice versa, ranks swap with degrees, etc. It is easy to see that $H$ is linear if and only if $H^*$ is linear. This is because $H$ is linear if and only if every $2 \times 2$ minor has a zero entry.

**Clique and line graphs**. The clique graph $C(H)$ has the same vertex set as $H$ and an edge for every pair of vertices in some edge. Each edge in $H$ then appears as a clique in $C(H)$. The line graph of $H$ is the clique graph of $H^*$; it has a vertex for each edge of $H$ and an edge between two edges of $H$ if they intersect. Note that in a linear hypergraph, the cliques in the clique graph which come from the edges in the hypergraph are edge disjoint; the clique graph is a set of edge disjoint complete subgraphs of the complete graph on $n$ vertices

**Coloring**. Suppose each edge of the hypergraph $H$ has a set of colors (a list) $L_e$ associated with it. A list edge coloring of $H$ is a function $\gamma(e)$ on the edges such that $\gamma(e) \in L_e$ and $\gamma(e) = \gamma(f)$ only $e$ and $f$ are independent. If $H$ can be list edge colored using any arbitrary set of lists as long as they each have $k$ colors, then we say that $H$ is $k$ list edge colorable. We let $q_{list}(H)$ be the smallest $k$ for which $H$ is $k$ list edge colorable. If all the lists are identical, we have an edge coloring and we call the smallest number of colors $q(H)$.

**Conjectures.** There is a long standing conjecture known now as the Erdős-Faber-Lovász conjecture which in its edge formulation says

*Conjecture (EFL). Let $H$ be a linear hypergraph with $n$ vertices and no rank 1 edges. Then $q(H) \leq n$.*

In [1], we discussed relationships between many generalizations of this conjecture. Here we extend these results to list edge coloring with essentially the same proofs.

*Conjecture C1 (list-EFL). Let $H$ be a linear hypergraph with $n$ vertices and no rank 1 edges. Then $q_{list}(H) \leq n$.*

*Definition.* We define the *clique degree* of a vertex $x$ by

$$D(x) = \sum_{e \supset x}(r(e) - 1).$$

This is the degree of $x$ in the clique graph. The dual concept is the *clique rank* of an edge

$$R(e) = \sum_{x \in e}(d(x) - 1).$$

*Conjecture C2. Let $H$ be a linear hypergraph with no rank 1 edges such that for every vertex $D(x) \leq k$. Then $q_{list}(H) \leq k + 1$.*

*Conjecture C3.  Let $H$ be a linear hypergraph with maximum rank $P$, maximum degree $\Delta$ and no rank 1 edges  then $q_{list}(H) \leq \Delta P - \max(\Delta, P) + 1$.*

*Conjecture C4 (weak Vizing conjecture).  For any graph $G$, if $d(x) \leq k$ then $q_{list}(G) \leq k + 1$.*

**Background Facts (Edge Formulations).**  In this section, we give two known theorems that we shall utilize to deduce relationships between the conjectures.  This requires giving the thoerem in a form dual to its usual presentation.

*Fact 1 (Greedy coloring - edge formulation).  For any hypergraph if $R(e) \leq R$ for every edge $e$ then $q_{list}(H) \leq R + 1$.*

*Proof.*  Color the edges in any order.  At each step, the edge which is about to be colored can only be adjacent to $R$ edges so at least one color must be unused in its associated list.

Note that this theorem also applies to any hypergraph whose edges can be ordered $e_i$ so that removing the edges $\{e_1, e_2, \cdots, e_i\}$ leaves a hypergraph such that $R(e_{i+1}) \leq R$.

Also note that Conjecture C2 is a natural extension of the weak Vizing conjecture.  If you think of $H$ as an edge disjoint union of cliques, then the hypothesis of Conjecture C2 is that the sum total of all the edges of all the cliques that meet at a fixed vertex $x$ is at most $k$.  The Vizing conjecture colors these edges so that they all have different colors.  Conjecture C2 demands that all the edges that belong to the same clique have the same color.

*Definition.*  Three edges $(e_1, e_2, e_3)$ in $H$ are called a *triangle with side $e_1$* if the edges all have pair-wise non-empty intersections.  We denote by $T(e)$ the number of triangles in $H$ with side $e$.

*Fact 2 (Vu Theorem [3] – edge formulation).  Let $H$ be a hypergraph with $R(e) \leq R$. Suppose there exists an $f > 1$ such that for every edge, $T(e) \leq R^2 / f$.  Then there exists a universal constant $c$ such that $q_{list}(H) \leq cR / \log f$.*

**Some Theorems**.  In this section we use the facts to prove theorems.

*Theorem 4.  C2 implies both C1 and C3.*

*Proof.*  Assume C2 is true.  To prove C3, let $H$ be a linear hypergraph with maximum rank $P$, maximum degree $\Delta$ and no rank 1 edges. We have two cases depending upon whether $P > \Delta$ or not. If $P > \Delta$, then we employ the Greedy coloring so that we have $q_{list}(H) \leq \Delta P - P + 1$.    If $P \leq \Delta$ then  $q_{list}(H) \leq \max D(x) \leq \Delta P - \Delta$.  This proves C3. To prove C1, let $H$ be a linear hypergraph with $n$ vertices and no rank 1 edges and let $x$ be a vertex.  Since $H$ is linear, all the vertices besides $x$ in all the edges that meet $x$ are distinct so there can be at most $n-1$ of them. Thus $D(x) \leq n-1$ and so $q(H) \leq n$. This proves C1.

*Theorem 5. Suppose that $H$ is a linear hypergraph of rank $\mathrm{P}$ with no rank 1 edges. There is a universal constant $C$ such that if $\mathrm{P} \geq C \geq 3$ and $\Delta \geq C(\mathrm{P}-1)$ then*

$$q_{list}(H) \leq \Delta(\mathrm{P}-1).$$

*Proof.* If $H$ is not uniform, add vertices of degree one as needed to make it uniform. Thus we can assume that $H$ is uniform with the rank of every edge equal to $r = \mathrm{P} \geq 3$. Let $k = \Delta(r-1)$ and assume that $H$ is not $k-1$ colorable. We work with the line graph $L = L(H)$ of the hypergraph $H$. The line graph uses the $m$ edges of $H$ as vertices and two edges in $H$ form an edge in $L$ if they meet in $H$. Any coloring of the vertices of $L$ is a coloring of the edges of $H$ and vice versa. We are given for every $x$

$$k \geq (r-1)d(x)$$

so

$$d(x) \leq \frac{k}{r-1}.$$

(1)

Thus

$$d(x) - 1 \leq \frac{k}{r-1} - 1 = \frac{k+1-r}{r-1}.$$

Thus degree $R(e)$ of a vertex $e$ in the line graph $L$ can be no more than

$$R(e) = \sum_{x \in e}(d(x)-1) \leq r\frac{k+1-r}{r-1}.$$

In particular, if we let

$$R = r\frac{k+1-r}{r-1} = (k+1)\frac{r}{r-1} - \frac{r^2}{r-1} \qquad (2)$$

the line graph has maximum degree at most $R$. Furthermore by Eq. 1 and Eq.2,

$$d(x) \leq \frac{R}{r} + 1.$$

(3)

By the greedy coloring, $L$ can be colored in at most $R+1$ colors. Thus we only have to deal with cases where $R \geq k$. So we assume that

(4) $$k \leq R.$$

We proceed by bounding the number $T(e)$ of triangles in $H$ with side $e$. This number is made up of two types of triangles which have edges $e$ and two other edges. First, there is $T_1(e)$, the triangles where all 3 edges have a common intersection. Second, there is $T_2(e)$, the triangles which meet pairwise at distinct vertices. The number of pairs of edges that meet a vertex $x$ in $e$ is given using Eq. 3 by

$$T_1(e) = \sum_{x \in e} \binom{d(x)-1}{2} \leq \frac{r}{2}\frac{R}{r}\left(\frac{R}{r}-1\right) = \frac{R}{2}\left(\frac{R}{r}-1\right).$$

To bound $T_2(e)$, we first select one of the $R(e)$ edges $e'$ that meets $e$ at some point $x$ and then select a pair of vertices $y \in e$ and $z \in e'$ to form the triangle with $x$. Since $H$ is linear these vertices can be in only one edge. This will count every triangle that exists twice so

$$T_2(e) \leq R\frac{(r-1)^2}{2}.$$

Thus

$$T(e) = T_1(e) + T_2(e) \leq \frac{R^2}{2}\left(\frac{1}{r} + \frac{(r-1)^2}{R}\right).$$

Let

$$\frac{1}{f} = \frac{1}{2r} + \frac{(r-1)^2}{2R}.$$

Then $T(e) \leq \frac{R^2}{f}$ so $q_{list}(H) \leq c\frac{R}{\log f}$ by Vu. Let $C = e^{2c}$. Since $r \geq C$, $2 - C/r \geq 1$.
Since $\Delta \geq C(r-1)$, utilizing inequality Equation 4 gives

$$R \geq k = \Delta(r-1) \geq C(r-1)^2 \geq \frac{C(r-1)^2}{2 - C/r}$$

or

$$\frac{(r-1)^2}{2R} \leq \frac{2 - C/r}{2C} = \frac{1}{C} - \frac{1}{2r}.$$

Thus

$$\frac{1}{f} = \frac{1}{2r} + \frac{(r-1)^2}{2R} \leq \frac{1}{C}.$$

Thus we have shown that under the given hypothesis that

$$q(H) \leq c\frac{R}{\log f} \leq \frac{R}{2}.$$

But from Equation 2

$$\frac{R}{2} = \frac{kr}{2(r-1)} + \frac{r}{2(r-1)} - \frac{r^2}{2(r-1)} \leq \frac{3}{4}k - \frac{r}{2} \leq \Delta(r-1).$$

*Corollary 6. Let $C$ be the universal constant from Theorem 5. Let $H$ be a hypergraph with $n$ vertices, $P \geq C \geq 3$, $\Delta \geq C(P-1)$ and no rank 1 edges. Then if $H$ is uniform or $n > \Delta(P-1)$, $q_{list}(H) \leq n-1$.*

*Proof.* For every $x$, $D(x) \leq n-1$. If $H$ is uniform,

$$d(x)(P-1) = D(x) \leq n-1$$

and so

$$n > \Delta(P-1).$$

*Theorem 7 (see [4]). If $n > (\Delta-1)^2$ or $n < \rho^2$ then $q_{list}(H) \leq n$.*

*Proof.* Suppose we start with a minimal (with respect to the number of edges) counterexample $H$. Let the maximum degree of $H$ be $\Delta \geq 2$. If an edge $e$ has rank less than $n/(\Delta-1)$, then by induction we color $H \setminus \{e\}$ using $n$ colors and then since the number of colors sharing a vertex with $e$ is less than $n$, we can extend the coloring to all of $H$. Thus we can assume $\rho \geq n/(\Delta-1)$. Now let $x$ be a vertex of degree $\Delta$. Count the vertices adjacent to $x$. The number is at least

$$\Delta\left(\frac{n}{\Delta-1} - 1\right) = n + \frac{n}{\Delta-1} - \Delta$$

but can be no more than $n-1$ so

$$n + \frac{n}{\Delta-1} - \Delta \leq n-1.$$

Thus $n \leq (\Delta-1)^2$. Similarly, by counting the number of vertices adjacent to any vertex $x$ of degree $d(x)$,

$$d(x)(\rho-1) \leq n-1$$

so

$$\Delta \leq \frac{n-1}{\rho-1}.$$

Thus

$$n \leq \left(\frac{n-\rho}{\rho-1}\right)^2.$$

Solve this quadratic to get

$$n \geq \rho^2.$$

**List edge coloring.**

*Theorem 8.* Suppose C4 holds. Let $H$ be a linear hypergraph with $n$ vertices. Let $H_3$ be the hypergraph with edges of rank less than 3 removed. Let $H_2$ be the hypergraph with only the edges of rank 2. If the maximum degree of $H_3$ is $\Delta$ and the maximum degree of $H_2$, $\Delta(H_2)$, is at most $n - 2\Delta - 1$ then any list $n$ coloring of the edges $H_3$ can be extended to a list $n$ coloring of $H$.

*Proof.* The number of colors that are not used by the edges at $x$ in the coloring of $H_3$ is at least $n - \Delta$. Given an edge $e = (x, y)$ in $H_2$, the number of colors available to color $e$ is at least

$$(n - \Delta) + (n - \Delta) - n = n - 2\Delta \leq \Delta(H_2) + 1.$$

If the weak Vizing list edge coloring conjecture holds, then each edge $e$ of $H_2$ can be colored using colors that are not used by other edges at the endpoints of $e$.

*Corollary 9.* Suppose C4 holds. If for every vertex $D(x, H_3) \geq 2\Delta$, then any $n$ list edge coloring of $H_3$ can be extended to an $n$ list edge coloring of $H$.

*Proof.* The degree of a vertex $x$ in the graph $H_2$ satisfies

$$d(x, H_2) = n - 1 - \sum_{e \supset x}(r(e) - 1) = n - 1 - D(x, H_3) \leq n - 1 - 2\Delta.$$

Thus $\Delta(H_2) \leq n - 1 - 2\Delta$.

*Corollary 10.* Suppose C4 holds. If $H_3$ is regular with degree $d$, then any $n$ list edge coloring of $H_3$ can be extended to an $n$ list edge coloring of $H$.

*Proof.* We can calculate $D(x, H_3) = \sum_{e \supset x}(r(e) - 1) \geq 2d$.

*Definition.* Let $d_k(x)$ be the number of edges with rank $k$ at the vertex $x$. Let $d(x)$ be the degree of the vertex $x$ in $H_3$. We call $\Delta - d(x)$ the *deficit* at $x$ and $\sum_{k \geq 4}(k - 3)d_k(x)$ the *excess* at $x$.

*Corollary 11.* Suppose C4 holds. If for every $x$ the excess at $x$ is at least as great as twice the deficit at $x$, then any $n$ list coloring of $H_3$ can be extended to an $n$ list coloring of $H$.

*Proof.* We have

$$d(x, H_2) = n - 1 - \sum_{e \supset x}(r(e) - 1) = n - 1 - 2d(x) - \sum_{k \geq 4}(k - 3)d_k(x) \leq n - 1 - 2\Delta.$$

**Critical hypergraphs**. Let $H$ be a linear hypergraph with maximum rank $P \geq 3$ and maximum degree $\Delta$. We let

$$D = D(H) = \max_x D(x) = \max_x \sum_{e \supset x}(r(e)-1).$$

We want to investigate what properties $H$ must have to be a minimal counterexample to C2.

*Theorem 12. If $H$ is a minimal counterexample to C2 then*

(i) $q_{list}(H) > D(H) + 1$;
(ii) for every $e$, $D(H \setminus e) = D(H)$;
(iii) for every $e$, $q_{list}(H \setminus e) = D(H) + 1$;
(iv) for every $e$, $R(e) \geq D$.

*Proof.* (ii) If $D(H \setminus e) < D(H)$ then the minimality of $H$ implies that $q_{list}(H \setminus e) \leq D(H)$. But then $q_{list}(H) \leq D + 1$.
(iii) Consequently, if $q_{list}(H \setminus e) - 1 < D(H) = D(H \setminus e)$ then $q_{list}(H) \leq D + 1$.
(iv) Suppose that $q_{list}(H) > D + 1$ but for every edge $e$, $q_{list}(H \setminus e) \leq D + 1$. If $e$ is an edge such that $R(e) < D$ then when the edges of $H \setminus e$ are colored in $D+1$ colors, only $D$ are used on the edges that meet $e$. Thus there is a color available for $e$ contradicting $q_{list}(H) > D + 1$.

**Summary**

We have discussed a number of ways to generalize the EFL conjecture to list coloring of arbitrary linear hypergraphs and exhibited various relationships between these new conjectures. In particular, for large enough fixed rank we can generalize EFL to obtain a generalization of Vizing's list color conjecture that we can show is true except for possibly a finite number of hypergraphs.